\documentclass[a4paper,12pt]{amsart}
\usepackage{amssymb}

\input xypic
\xyoption{all}
\xyoption{arc}

\newtheorem{theo}{Theorem}[section]

\newtheorem{prop}[theo]{Proposition}

\newtheorem{lem}[theo]{Lemma}

\newtheorem{coro}[theo]{Corollary}

\def\remark#1{{\refstepcounter{theo}\label{#1}\noindent\sc Remark
\arabic{section}.\arabic{theo} - }}

\def\equat{\refstepcounter{theo}$$~}
\def\endequat{\leqno{\boldsymbol{(\arabic{section}.\arabic{theo})}}~$$}

    \def\FM{{\mathbb{F}}}

  \def\lG{{\mathfrak l}}

    \def\QM{{\mathbb{Q}}}

    \def\ZM{{\mathbb{Z}}}

    \def\AC{{\mathcal{A}}}
\def\Bb{{\mathbf B}}    \def\BC{{\mathcal{B}}}
    
\def\Db{{\mathbf D}}    
    \def\EC{{\mathcal{E}}}
    
\def\Gb{{\mathbf G}}    \def\GC{{\mathcal{G}}}
    \def\HC{{\mathcal{H}}}

\def\Lb{{\mathbf L}}    
\def\Mb{{\mathbf M}}    
    
    \def\OC{{\mathcal{O}}}
\def\Pb{{\mathbf P}}    
    
    \def\RC{{\mathcal{R}}}
\def\Sb{{\mathbf S}}    \def\SC{{\mathcal{S}}}
\def\Tb{{\mathbf T}}    
\def\Ub{{\mathbf U}}    
\def\Vb{{\mathbf V}}

\def\Yb{{\mathbf Y}}    
\def\Zb{{\mathbf Z}}    \def\ZC{{\mathcal{Z}}}

\def\Rrm{{\mathrm{R}}}    
    
    \def\TCB{{\boldsymbol{\mathcal{T}}}}

\def\Ati{{\tilde{A}}}

  \def\hti{{\tilde{h}}}

  \def\sti{{\tilde{s}}}  
  \def\tti{{\tilde{t}}}

  \def\zti{{\tilde{z}}}  

\def\Bbt{{\tilde{\Bb}}}

\def\Gbt{{\tilde{\Gb}}}

\def\Lbt{{\tilde{\Lb}}}

\def\Pbt{{\tilde{\Pb}}}

\def\Tbt{{\tilde{\Tb}}}

\def\a{\alpha}

\def\g{\gamma}
\def\G{\Gamma}
\def\d{\delta}

\def\e{\varepsilon}

\def\L{\Lambda}

\def\s{\sigma}
\def\th{\theta}

\def\t{\tau}

        \def\alpt{{\tilde{\alpha}}}

\def\Gamb{{\boldsymbol{\Gamma}}}

          \def\taut{{\tilde{\t}}}

\DeclareMathOperator{\Hom}{{\mathrm{Hom}}}

\DeclareMathOperator{\Id}{{\mathrm{Id}}}

\DeclareMathOperator{\im}{{\mathrm{Im}}}

\DeclareMathOperator{\Ind}{{\mathrm{Ind}}}

\DeclareMathOperator{\Irr}{{\mathrm{Irr}}}

\DeclareMathOperator{\Tr}{{\mathrm{Tr}}}

\DeclareMathOperator{\Cur}{Cur}
\DeclareMathOperator{\KCur}{\lind{\mathit{K}}{Cur}}
\DeclareMathOperator{\CurSh}{\Delta}
\DeclareMathOperator{\KCurSh}{\lind{\mathit{K}}{\Delta}}
\DeclareMathOperator{\Norm}{N}
\DeclareMathOperator{\KNorm}{\lind{\mathit{K}}{N}}
\DeclareMathOperator{\KInd}{\lind{\mathit{K}}{Ind}}

\def\to{\rightarrow}
\def\longto{\longrightarrow}
\def\injto{\hookrightarrow}

\def\longmapright#1{\hspace{0.3em}\smash{
     \mathop{\longrightarrow}\limits^{#1}}\hspace{0.3em}}

\def\fonction#1#2#3#4#5{\begin{array}{rccc}
{#1} : & {#2} & \longto & {#3} \\
& {#4} & \longmapsto & {#5}
\end{array}}

\def\DS{\displaystyle}
\def\SS{\scriptstyle}

\def\finl{~$\SS \square$}

\def\lexp#1#2{\kern\scriptspace\vphantom{#2}^{#1}\kern-\scriptspace#2}
\def\lind#1#2{\kern\scriptspace\vphantom{#2}_{#1}\kern-\scriptspace#2}

\def\ge{\hspace{0.1em}\mathop{\geqslant}\nolimits\hspace{0.1em}}

\mathchardef\inferieur="321E
\mathchardef\superieur="321F

\def\eqna{\begin{eqnarray*}}
\def\endeqna{\end{eqnarray*}}

\def\itemth#1{\item[${\mathrm{(#1)}}$]}

\def\gfp{\FM_{\! p}}
\def\ql{\QM_{\!\ell}}
\def\sem{{\mathrm{sem}}}
\def\gfq{\FM_{\! q}}

\begin{document}

\baselineskip=16pt

\title{On the endomorphism algebras of  modular Gelfand-Graev representations}

\author{C\'edric Bonnaf\'e$^1$}
\address{$^1$ Labo. de Math. de Besan\c{c}on
(CNRS: UMR 6623),
Universit\'e de Franche-Comt\'e, 16 Route de Gray, 25030 Besan\c{c}on
Cedex, France}

\makeatletter
\email{cedric.bonnafe@univ-fcomte.fr}

\author{Radha Kessar$^2$}
\address{$^2$ Mathematical Sciences, University of Aberdeen,
Meston Building, Aberdeen AB24 3UE, Scotland, UK}

\makeatletter
\email{kessar@maths.abdn.ac.uk}
\makeatother


\date{\today}
\maketitle

\centerline{\it This paper is dedicated to Toshiaki Shoji, on his sixtieth birthday.}

\begin{abstract}
We study the endomorphism algebras of a modular Gelfand-Graev representation
of a finite reductive group by investigating modular properties of
homomorphisms constructed by Curtis and Curtis-Shoji.
\end{abstract}

\pagestyle{myheadings}

\markboth{\sc C. Bonnaf\'e \& R. Kessar}{\sc Endomorphism algebras of
Gelfand-Graev representations}

\vskip1cm



Let $\Gb$ be  a connected reductive group defined over an algebraic closure
$\FM$ of the  field of $p$-elements $\gfp$
and suppose that it is endowed with a Frobenius endomorphism  $F : \Gb \to \Gb$
relative to  an  $\gfq$-structure.  Since   the work of  Lusztig,
it  has been  natural   to ask   to what extent the   theory  of the
representations  of  $\Gb^F$  depends on   $q$.  For example, it was shown by Lusztig  that
the unipotent characters of  $\Gb^F$  are parametrized by a set  which is independent of  $q$  (the set depends solely on the  Weyl group of  $\Gb$  and on the action of  $F$  on this  Weyl  group).

On the side of  $\ell$-modular  representations  (where  $\ell$ is a prime different from   $p$),  the work  of  Fong and Srinivasan
on the  general linear and unitary groups \cite{fong1} and
on the classical groups \cite{fong}, then that of
Brou\'e, Malle and Michel (introducing the notion of  {\it generic groups} \cite{BMM}) and of Cabanes and Enguehard \cite{cabanes1}  
give evidence of analogous results.  For instance,  in  most cases,
the  unipotent  $\ell$-blocks of $\Gb^F$  only depend on the order of
$q$ modulo $\ell$, and not on the value of  $q$ itself \cite[Chapter 22]{cabanes}.

Let  $(K,\OC,k)$  denote an  $\ell$-modular  system,  sufficiently
large. In this article we  will study the endomorphism
algebra $\HC^\Gb_{(d)}$  of a modular  Gelfand-Graev representation
 $\G^\Gb_{(d)}$ of $\Gb^{F^d}$ (this is a projective  $\OC\Gb^{F^d}$-module).
We will study also  the corresponding   unipotent parts
$b_{(d)}^\Gb \HC^\Gb_{(d)}$ and $b_{(d)}^\Gb \G^\Gb_{(d)}$
(here, $b_{(d)}^\Gb$ denotes  the sum of the unipotent blocks  of $\OC\Gb^{F^d}$).  We make the following conjecture, which is   related
to the question   mentioned above:

\bigskip

\begin{quotation}
\noindent{\bf Conjecture 1.}
{\it  If $\ell$ does not  divide $[\Gb^{F^d}:\Gb^F]$,  then the
$\OC$-algebras
$b_{(d)}^\Gb \HC^\Gb_{(d)}$ and  $b_{(1)}^\Gb \HC^\Gb_{(1)}$  are
isomorphic.}
\end{quotation}

\bigskip

The  $\OC\Gb^{F^d}$-module $b_{(d)}^\Gb \G_{(d)}^\Gb$  is
projective and  indecomposable~: it is the  projective cover  of the  modular {\it  Steinberg  module}.
Conjecture 1,  if proven, would show that
the endomorphism algebra  of this   module  does not depend
{\it too much } on $q$.

In this article,  we approach Conjecture  1 by the study of a morphism
$\KCur_{\Lb,(d)}^\Gb : K\HC^\Gb_{(d)} \to K\HC^\Lb_{(d)}$
(where  $\Lb$  is an $F^d$-stable  Levi subgroup  of a parabolic subgroup of
 $\Gb$).  When  $\Tb$  is a maximal
$F^d$-stable  torus of $\Gb$,
 this morphism   was constructed by  Curtis \cite[Theorem 4.2]{curtis}
and it is defined over  $\OC$ (i.e.  there exists a morphism of algebras
$\Cur_{\Tb,(d)}^\Gb : \HC^\Gb_{(d)} \to \HC^\Tb_{(d)}=\OC\Tb^{F^d}$
such  that $\KCur_{\Tb,(d)}^\Gb$ is  obtained  from $\Cur_{\Tb,(d)}^\Gb$ by extension of scalars).
We will  also  consider a product of
 Curtis homomorphisms
$$\Cur_{(d)}^\Gb : \HC_{(d)}^\Gb \longto \prod_{\Tb \in \TCB^{F^d}}
\OC\Tb^{F^d},$$
where  $\TCB$  is the variety of maximal tori of  $\Gb$.
Finally,   we   will study  a morphism of $K$-algebras
$\KCurSh^\Gb : K\HC^\Gb_{(d)} \to K\HC^\Gb_{(1)}$
defined by  Curtis and Shoji \cite[Theorem 1]{curtis shoji}.
With this notation, we can  state Conjecture 1 more precisely~:

\bigskip

\begin{quotation}
\noindent{\bf Conjecture 2.}
{\it  With  the notation above, we have:
\begin{itemize}
\itemth{a} $\KCur_\Lb^\Gb$  is defined over $\OC$.

\itemth{b} $\KCurSh^\Gb$ is defined  over $\OC$.

\itemth{c} If  $\ell$  does not divide $[\Gb^{F^d}:\Gb^F]$, then
$\KCurSh^\Gb$ induces an isomorphism
$b_{(d)}^\Gb \HC^\Gb_{(d)} \simeq b_{(1)}^\Gb \HC^\Gb_{(1)}$.
\end{itemize}}
\end{quotation}

\bigskip

The main results of this article are obtained  under the hypothesis  that
$\ell$  does not divide the order of the
 Weyl group $W$ of $\Gb$.

\bigskip

\noindent{\bf Theorem.}
{\it If  $\ell$  does not divide  $|W|$,  then  Conjecture 2 holds.}

\bigskip

Statement (a)  is  proved in Corollary \ref{cur O}~;
statement (b)   in   Theorem \ref{curtis shoji over dvr}~;
statement (c) is shown in  Theorem \ref{endo iso}.
In order to obtain our theorem,  we proved two   more precise  results
which   do not necessarily  hold   when $\ell$ does divide $|W|$.

\bigskip

\noindent{\bf Theorem \ref{image}.}
{\it If  $\ell$ does not divide $|W|$, then
$$\im(\Cur_{(d)}^\Gb)=\im(\KCur_{(d)}^\Gb)
\cap \prod_{\Tb \in \TCB^{F^d}} \OC\Tb^{F^d}.$$}

\bigskip

\noindent{\bf Theorem \ref{image sylow}.}
{\it If  $\ell$  does not divide  $|W|$, then
$$b_{(d)}^\Gb \HC_{(d)}^\Gb \simeq (\OC S)^{N_{\Gb^{F^d}}(S)},$$
where  $S$ is a Sylow $\ell$-subgroup of $\Gb^{F^d}$.}

\noindent{\sc Remark-} With the above notation, if $\ell$  does not divide  $|W|$, then 
$S$ is abelian, and hence a     consequence of  the above  result 
is that   if $\ell$  does not divide  $|W|$, then  
the isomorphism type of the $\OC$-algebra 
$b_{(d)}^\Gb \HC_{(d)}^\Gb$  depends only on the fusion of $\ell $-elements
in $\Gb^{F^d}$.

\bigskip

This article is organized as follows.
In the first section, we  recall the definitions of the
Gelfand-Graev representations  as well  as some of the principal
properties of their  endomorphism algebras  (commutativity for example).
In the second section,  we construct the generalisation of the Curtis
homomorphism. In the third part we study the product of  Curtis
homomorphisms and  prove, amongst other things,  Theorems \ref{image}
and  \ref{image sylow} stated above.
In the last part, we study the Curtis-Shoji homomorphism and prove
statement (c)  of Conjecture 2 when $ \ell$  does not divide $|W|$.

\vskip1cm

{

\noindent{\sc Notation - } If $A$ is a finite dimensional algebra
over a field, we denote by $\RC(A)$ the Grothendieck group of the
category of finitely generated $A$-modules. If $M$ is a finitely
generated $A$-module, we denote by $[M]$ its class in $\RC(A)$.
The opposite algebra of $A$ will be denoted by $A^\circ$.

All along this paper, we fix a prime number $p$, an algebraic closure
$\FM$ of the finite field with $p$ elements $\gfp$, a prime number
$\ell$ different from $p$ and an algebraic extension $K$ of the $\ell$-adic field
$\ql$. Let $\OC$ be the ring of integers of $K$, let $\lG$ be the
maximal ideal of $\OC$ and let $k$ denote the residue field of
$\OC$: $k$ is an algebraic extension of the finite field $\FM_{\! \ell}$.
Throughout this paper, we assume that the $\ell$-modular system
$(K,\OC,k)$ is sufficiently large for all the finite groups
considered in this paper.

If $\L$ is a commutative $\OC$-algebra (for instance $\L=k$ or $K$), and if $M$
is an $\OC$-module, we set $\L M=\L \otimes_\OC M$. If $f : M \to N$
is a morphism of $\OC$-modules, we define $\lind{\L}{f} : \L M \to \L N$
to be the morphism $\Id_\L \otimes_\OC f$. If $V$ is a free left $\L$-module,
we denote by $V^*=\Hom_\L(V,\L)$ its dual: if $V$ is a left $A$-module for some
$\L$-algebra $A$, then $V^*$ is seen as a right $A$-module.

If $G$ is a finite group, we denote by $\Irr G$ the set of irreducible
characters of $G$ over $K$. If $\chi \in \Irr G$, let
$e_\chi$ (or $e_\chi^G$ if we need to emphasize the ambient group)
denote the associated central primitive idempotent of $KG$~:
$$e_\chi=\frac{\chi(1)}{|G|}\sum_{g \in G} \chi(g^{-1}) g.$$
The conjugacy relation in $G$ is denoted by $\sim$ or $\sim_G$ if
necessary.

\vskip1cm

\section{Background material}

\medskip

\subsection{The set-up}
We fix once and for all a connected reductive algebraic group
$\Gb$ over $\FM$ and we assume that it is endowed with an isogeny $F : \Gb \to \Gb$
such that some power of $F$ is a Frobenius endomorphism of $\Gb$
with respect to some rational structure on $\Gb$ over a finite extension
of $\gfp$. We denote by $q$ the positive real number such that,
for every $\d \ge 1$ such that $F^\d$ is a Frobenius endomorphism
of $\Gb$ over a finite field with $r$ elements, we have $r=q^\d$.

\bigskip

\subsection{Gelfand-Graev representations}
We fix an $F$-stable Borel subgroup $\Bb_\Gb$ of $\Gb$ and an $F$-stable
maximal torus of $\Bb_\Gb$. Let $\Ub_\Gb$ denote the unipotent radical
of $\Bb_\Gb$. We fix once and for all a {\it regular} linear character
$\psi : \Ub_\Gb^F \to \OC^\times \subset K^\times$ (in the sense of
\cite[Definition 2.3]{DLM1}).
Since $|\Ub_\Gb^F|=q^{\dim \Ub_\Gb}$
is a power of $p$, the primitive central idempotent $e_\psi$
of $K\Ub_\Gb^F$ belongs to $\OC\Ub_\Gb^F$.
We denote by $\OC_\psi$ the projective
$\OC\Ub_\Gb^F$-module $\OC\Ub_\Gb^Fe_\psi$: it is $\OC$-free of rank one and
is acted on by $\Ub_\Gb^F$ through $\psi$. Let
$$\G^\Gb = \OC\Gb^F e_\psi \simeq \Ind_{\Ub_\Gb^F}^{\Gb^F} \OC_\psi.$$
Then $\G^\Gb$ is a projective $\OC\Gb^F$-module; the corresponding
representation is called a {\it Gelfand-Graev representation} of $\Gb^F$.

Let $\HC^\Gb$ denote the endomorphism algebra of the $\OC\Gb^F$-module
$\G^\Gb$. We have
\equat\label{endo gelfand}
\HC^\Gb \simeq (e_\psi \OC\Gb^F e_\psi)^\circ
\endequat
Since $\OC\Gb^F$ is a symmetric algebra, we have that
\equat\label{h symetrique}
\text{\it $\HC^\Gb$ is symmetric.}
\endequat
The next result is much more difficult (see \cite[Theorem 15]{steinberg} for the
general case):

\begin{theo}
The algebra $\HC^\Gb$ is commutative.
\end{theo}

Therefore, we shall identify the algebras $\HC^\Gb$ and $e_\psi \OC\Gb^F e_\psi$.

\medskip

If $\L$ is a commutative $\OC$-algebra, we denote by $e_\psi^\L$ the idempotent
$1_\L \otimes_\OC e_\psi$ of $\L \Ub_\Gb^F=\L \otimes_\OC \OC\Ub_\Gb^F$.
Since $\G^\Gb$ is projective, the $\L\Gb^F$-module $\L\G^\Gb$ is also
projective and its endomorphism algebra is $\L\HC^\Gb$ (since
it is isomorphic to $\Hom_\L(\L\G^\Gb,\L) \otimes_{\L\Gb^F} \L \G^\Gb$).
We have of course (taking into account that $\HC^\Gb$ is symmetric)
\equat\label{endo gelfand bis}
\L\HC^\Gb = e_\psi^\L \L\Gb^F e_\psi^\L.
\endequat
Since $K\Gb^F$ is split semisimple,
\equat\label{split}
\text{\it The algebra $K\HC^\Gb$ is split semisimple.}
\endequat


\remark{several}
There might be several Gelfand-Graev representations of $\OC\Gb^F$. But they
are all conjugate by elements $g \in \Gb$ such $g^{-1}F(g)$ belongs to the
centre of $\Gb$, and this gives a parametrization of Gelfand-Graev
representations by the group of $F$-conjugacy classes in the centre of $\Gb$
(see \cite[2.4.10]{DLM1}). In particular, their endomorphism algebras
are all isomorphic.

Moreover, if the centre of $\Gb$ is connected, there
is only one (up to isomorphism) Gelfand-Graev representation.
In special orthogonal or symplectic groups in odd characteristic, there
are two (isomorphism  classes of) Gelfand-Graev representations. \finl

\bigskip

\subsection{Representations of ${\boldsymbol{K\HC^\Gb}}$}
Let $(\Gb^*,F^*)$ be a dual pair to $(\Gb,F)$ in the sense of
Deligne and Lusztig \cite[Definition 5.21]{delu}. We denote by $\Gb^*_\sem$ the
set of semisimple elements of $\Gb^*$. If $s \in \Gb_\sem^{*F^*}$, we
denote by $(s)_{\Gb^{*F^*}}$ its conjugacy class in $\Gb^{*F^*}$ and
by $\EC(\Gb^F,(s)_{\Gb^{*F^*}})$ the associated rational Lusztig
series (see \cite[Page 136]{dmbook}). We denote by $\chi_s^\Gb$ the unique
element of $\EC(\Gb^F,(s)_{\Gb^{*F^*}})$ which is an irreducible
component of the character afforded by $K\G^\Gb$. We view it as
a function $K\Gb^F \to K$ and we denote
by $\chi_s^{\HC^\Gb}$ its restriction to $K\HC^\Gb$.
Then (see \cite[Theorem 10.7]{delu} for the case where the centre
of $\Gb$ is connected and \cite{asai} for the general case; see also
\cite[Remark of Page 80]{bonnafe asterisque} for the
case where $F$ is not a Frobenius endomorphism),
\equat\label{decomposition}
[K\G^\Gb] = \sum_{(s)_{\Gb^{*F^*}} \in \Gb_\sem^{*F^*}/\sim} \chi_s^\Gb.
\endequat
Therefore, the next proposition is a particular case of \cite[Theorem 11.25
and Corollaries 11.26 and 11.27]{CR},
taking into account that $K\HC^\Gb$ is semisimple and commutative
(or, equivalently, that $K\G^\Gb$ is multiplicity free):

\begin{prop}\label{representations H}
We have:
\begin{itemize}
\itemth{a} The map $s \mapsto \chi_s^{\HC^\Gb}$ induces a bijection between
the set of $\Gb^{*F^*}$-conjugacy classes of semisimple elements of
$\Gb^{*F^*}$ and the set of irreducible characters of $K\HC^\Gb$.

\itemth{b} The map $s \mapsto e_{\chi_s^\Gb}e_\psi$ induces a bijection between
the set of $\Gb^{*F^*}$-conjugacy classes of semisimple elements of
$\Gb^{*F^*}$ and the set of primitive idempotents of $K\HC^\Gb$.
\end{itemize}
Moreover, if $s \in \Gb_\sem^{*F^*}$, then:
\begin{itemize}
\itemth{c} We have $\chi_s^{\HC^\Gb}(e_{\chi_s^\Gb} e_\psi)=1$
and $\chi_s^{\HC^\Gb}(e_{\chi_t^\Gb} e_\psi)=0$ if $t \in \Gb_\sem^{*F^*}$
is not conjugate to $s$ in $\Gb^{*F^*}$.

\itemth{d} $e_{\chi_s^\Gb} e_\psi$ is the primitive idempotent of $K\HC^\Gb$
associated with the irreducible character $\chi_s^{\HC^\Gb}$.

\itemth{e} The $K\Gb^F$-module
$K\Gb^Fe_{\chi_s^\Gb} e_\psi$ is irreducible and affords the character $\chi_s^\Gb$.

\itemth{f} If $\chi \in \EC(\Gb^F,(s)_{\Gb^{*F^*}})$ and if $\chi \neq \chi_s^\Gb$,
then $\chi(h)=0$ for all $h \in K\HC^\Gb$.
\end{itemize}
\end{prop}

Since $K\HC^\Gb$ is split and commutative, all its
irreducible representations have dimension one. In other words,
all its irreducible characters are morphisms of $K$-algebras $K\HC^\Gb \to K$.
So, as a consequence of the Proposition \ref{representations H}, we get that
the map
\equat\label{injectif}
\fonction{\chi^{\HC^\Gb}}{K\HC^\Gb}{\displaystyle{\prod_{(s)_{\Gb^{*F^*}} \in
\Gb_\sem^{*F^*}/\sim} K}}{h}{(\chi_s^{\HC^\Gb}(h))_{(s)_{\Gb^{*F^*}} \in
\Gb_\sem^{*F^*}/\sim}}
\endequat
is an isomorphism of $K$-algebras. It corresponds to the decomposition
\equat\label{base idem}
K\HC^\Gb=\mathop{\oplus}_{(s)_{\Gb^{*F^*}} \in \Gb_\sem^{*F^*}/\sim}
K\HC^\Gb e_{\chi_s^\Gb} e_\psi.
\endequat

\bigskip

\section{A generalization of  the Curtis homomorphisms}

\medskip

 In   \cite[Theorem 4.2]{curtis}, Curtis
 constructed a homomorphism of algebras $f_\Tb : \HC^\Gb \to \OC\Tb^F$, for
$\Tb$  an $F$-stable maximal torus of $\Gb$
(in fact, Curtis  constructed a homomorphism of algebras $K\HC^\Gb \to K\Tb^F$
but it is readily checked from his formulas that it is defined
over $\OC$). We propose here a generalization of this construction to the case
where $\Tb$ is replaced by an $F$-stable Levi subgroup $\Lb$ of a parabolic subgroup
of $\Gb$: we then get a morphism $K\HC^\Gb \to K\HC^\Lb$
(note that, if $\Lb$ is a maximal torus, then $\HC^\Lb=\OC\Tb^F$).
We conjecture that this morphism is defined over $\OC$ and prove it
whenever $\GC^+(\Gb,\Lb,\Pb)$ holds or whenever $\Lb$ is a maximal
torus (see Theorem \ref{main curtis}) or whenever $\ell$ does not divide
the order of $W$ (see Corollary \ref{cur O}).

\bigskip

\subsection{A morphism ${\boldsymbol{K\HC^\Gb \to K\HC^\Lb}}$}
Let $\Pb$ be a parabolic subgroup of $\Gb$ and assume that $\Pb$ admits an $F$-stable
Levi complement $\Lb$.

The Gelfand-Graev representation $\G^\Gb$ of $\OC\Gb^F$ having been fixed,
there is a well-defined (up to isomorphism) Gelfand-Graev representation $\G^\Lb$ of
$\OC\Lb^F$  associated to it  \cite[Page 77]{bonnafe asterisque}
(see also \cite{bonnafe regular}).
We fix an $F$-stable Borel subgroup $\Bb_\Lb$ of $\Lb$ and we denote by
$\Ub_\Lb$ its unipotent radical. We fix once and for all a regular linear
character $\psi_\Lb$ of $\Ub_\Lb^F$ such that
$\G^\Lb=\Ind_{\Ub_\Lb^F}^{\Lb^F} \OC_{\psi_\Lb} = \OC\Lb^F e_{\psi_\Lb}$.
We identify $\HC^\Lb$ with $e_{\psi_\Lb} \OC\Lb^F e_{\psi_\Lb}$.
We also fix an $F^*$-stable Levi subgroup $\Lb^*$ of a parabolic subgroup of $\Gb^*$
dual to $\Lb$ (this  is well-defined up to conjugacy by an element of
$\Gb^{*F^*}$: see \cite[Page 113]{dmbook}). We then define
$\KCur_\Lb^\Gb K \HC^\Gb \to K\HC^\Lb$ as the unique linear map such that, for
any semisimple element $s \in \Gb^{*F^*}$,
$$\KCur_\Lb^\Gb(e_{\chi_s^\Gb} e_\psi) =
\sum_{\substack{(t)_{\Lb^{*F^*}} \in \Lb_\sem^{*F^*}/\sim \\ t \in (s)_{\Gb^{*F^*}}}}
e_{\chi_t^\Lb} e_{\psi_\Lb}.$$
Note that this does not depend on the choice of the representative $s$
in its conjugacy class.

\begin{prop}\label{prop curtis}
The map $\KCur_\Lb^\Gb$ is an homomorphism of algebras. Moreover, if
$s \in \Lb_\sem^{*F^*}$, then
$$\chi_s^{\HC^\Lb} \circ \KCur_\Lb^\Gb = \chi_s^{\HC^\Gb}.$$
\end{prop}

\begin{proof}
Since the image of an idempotent is an idempotent
(and since $K\HC^\Gb$ and $K\HC^\Lb$ are split semisimple and commutative),
we get the first statement. The second  is obtained by
applying both sides to each primitive idempotent $e_{\chi_t^\Gb} e_\psi$ of $K\HC^\Gb$
($t \in \Gb_\sem^{*F^*}$).
\end{proof}

Another easy consequence of the definition is the following

\begin{prop}\label{transitivity}
If $\Mb$ is an $F$-stable Levi subgroup of a
parabolic subgroup of $\Gb$ and if $\Lb \subset \Mb$, then
$\KCur_\Lb^\Mb \circ \KCur_\Mb^\Gb = \KCur_\Lb^\Gb$.
\end{prop}

\bigskip

\subsection{Deligne-Lusztig functors and Gelfand-Graev
representations\label{subsection dl}}
Let $\Pb$ be a parabolic subgroup of $\Gb$ and assume that $\Pb$ admits an $F$-stable
Levi complement $\Lb$. Let $\Vb$ denote the unipotent radical of $\Pb$. We set
$$\Yb_\Pb^\Gb =\{g \Vb \in \Gb/\Vb~|~g^{-1}F(g) \in \Vb \cdot F(\Vb)\}$$
and $d_\Pb=\dim(\Vb) - \dim(\Vb \cap F(\Vb))$. Then
$\Yb_\Pb^\Gb$ is a locally closed smooth variety of pure dimension
$d_\Pb$. If $\L=\OC$, $K$ or $\OC/\lG^n$,
the complex of cohomology with compact support of $\Yb_\Pb^\Gb$ with
coefficients in $\L$, which is denoted by $\Rrm\G_c(\Yb_\Pb^\Gb,\L)$,
is a bounded complex of $(\L\Gb^F,\L\Lb^F)$-bimodules which is perfect as a complex of
left $\L\Gb^F$-modules and is also perfect as a complex of right $\L\Lb^F$-modules
(see \cite[\S 3.8]{delu}).
Its $i$-th cohomology group is denoted by $H_c^i(\Yb_\Pb^\Gb,\L)$:
it is a $(\L\Gb^F,\L\Lb^F)$-bimodule. For $\L=\OC$, this
complex induces two functors between bounded derived categories
$$\fonction{\RC_{\Lb \subset \Pb}^\Gb}{D^b(\OC\Lb^F)}{D^b(\OC\Gb^F)}{C}{
\Rrm\G_c(\Yb_\Pb^\Gb,\OC) \otimes_{\OC\Lb^F} C}$$
$$\fonction{\lexp{*}{\RC}_{\Lb \subset \Pb}^\Gb}{D^b(\OC\Gb^F)}{D^b(\OC\Lb^F)}{C}{
\Rrm\Hom_{\OC\Gb^F}^\bullet(\Rrm\G_c(\Yb_\Pb^\Gb,\OC), C).}\leqno{\text{and}}$$
These functors are respectively called {\it Deligne-Lusztig induction}
and {\it restriction}.
By extending the scalars to $K$, they induce linear maps between
the Grothendieck groups $R_{\Lb \subset \Pb}^\Gb : \RC(K\Lb^F) \to \RC(K\Gb^F)$ and
$\lexp{*}{R}_{\Lb \subset \Pb}^\Gb : \RC(K\Gb^F) \to \RC(K\Lb^F)$.
We have
$$R_{\Lb \subset \Pb}^\Gb[M]= \sum_{i \ge 0} (-1)^i
[H^i_c(\Yb_\Pb^\Gb,K) \otimes_{K\Lb^F} M]$$
$$\lexp{*}{R}_{\Lb \subset \Pb}^\Gb[N] =\sum_{i \ge 0} (-1)^i
[H^i_c(\Yb_\Pb^\Gb,K)^* \otimes_{K\Gb^F} N]\leqno{\text{and}}$$
for all $K\Gb^F$-modules $N$ and all $K\Lb^F$-modules $M$.

If $(g,l) \in K\Gb^F \times K\Lb^F$, we set
$$\Tr_{\Lb \subset \Pb}^\Gb(g,l)=\sum_{i \ge 0} (-1)^i \Tr((g,l),H^i_c(\Yb_\Pb^\Gb,K)).$$
If $(g,l) \in \Gb^F \times \Lb^F$, then
$\Tr_{\Lb \subset \Pb}^\Gb(g,l)$ is a rational integer which does not
depend on the prime number $\ell$ (see \cite[Proposition 3.3]{delu}).
If $\chi_M$ (respectively $\chi_N$)
denotes the character afforded by a $K\Lb^F$-module $M$ (respectively a
$K\Gb^F$-module $N$), then the character afforded by the virtual module
$R_{\Lb \subset \Pb}^\Gb[M]$ (respectively $\lexp{*}{R}_{\Lb \subset \Pb}^\Gb[N]$)
will be denoted by $R_{\Lb \subset \Pb}^\Gb \chi_M$ (respectively
$\lexp{*}{R}_{\Lb \subset \Pb}^\Gb \chi_N$): it satisfies
$$R_{\Lb \subset \Pb}^\Gb \chi_M(g) = \frac{1}{|\Lb^F|} \sum_{l \in \Lb^F}
\Tr_{\Lb \subset \Pb}^\Gb(g,l) \chi_M(l^{-1})$$
$$\lexp{*}{R}_{\Lb \subset \Pb}^\Gb \chi_N(l) = \frac{1}{|\Gb^F|}\sum_{g \in \Gb^F}
\Tr_{\Lb \subset \Pb}^\Gb(g,l) \chi_N(g^{-1})\quad)\leqno{\text{(respectively}}$$
for all $g \in \Gb^F$ (respectively $l \in \Lb^F$).

\bigskip

\noindent{\sc Comments (independence on the parabolic) - }
If $\Pb'$ is another parabolic subgroup of $\Gb$ having $\Lb$ as a Levi
complement, then the Deligne-Lusztig varieties $\Yb_\Pb^\Gb$ and
$\Yb_{\Pb'}^\Gb$ are in general non-isomorphic: they might even have different
dimension (however, note that $(-1)^{d_\Pb}=(-1)^{d_{\Pb'}}$, i.e.
$d_\Pb \equiv d_{\Pb'} \mod 2$). As a consequence, the
Deligne-Lusztig functors $\RC_{\Lb \subset \Pb}^\Gb$
and $\RC_{\Lb \subset \Pb'}^\Gb$ can be really different.
However, it is conjectured in general that
$R_{\Lb \subset \Pb}^\Gb = R_{\Lb \subset \Pb'}^\Gb$ and
$\lexp{*}{R}_{\Lb \subset \Pb}^\Gb = \lexp{*}{R}_{\Lb \subset \Pb'}^\Gb$.
This is equivalent to say that $\Tr_{\Lb \subset \Pb}^\Gb =\Tr_{\Lb \subset \Pb'}^\Gb$.

For instance, we have $\Tr_{\Lb \subset \Pb}^\Gb =\Tr_{\Lb \subset \Pb'}^\Gb$
if $\Lb$ is a maximal torus \cite[Corollary 4.3]{delu}, or if $\Pb$ and $\Pb'$
are $F$-stable (this is due to Deligne: a proof can be found in
\cite[Theorem 5.1]{dmbook}), or if $F$
is a Frobenius endomorphism and $q \neq 2$ (see \cite{bonnafe mackey}).
In all these cases, this fact is a consequence of the Mackey formula for
Deligne-Lusztig maps.\finl

\bigskip

The Gelfand-Graev representation $\G^\Lb$
satisfies the following property:

\begin{theo}\label{rappel restriction}
Assume that one of the following three conditions is satisfied:
\begin{itemize}
\itemth{1} $\Pb$ is $F$-stable.

\itemth{2} The centre of $\Lb$ is connected.

\itemth{3} $p$ is almost good for $\Gb$, $F$ is a Frobenius endomorphism of
$\Gb$ and $q$ is large enough.
\end{itemize}
Then $\lexp{*}{R}_{\Lb \subset \Pb}^\Gb [K\G^\Gb]=(-1)^{d_\Pb} [K\G^\Lb]$.
\end{theo}

\begin{proof}
(1) is due to Rodier: a proof
may be found in \cite[Theorem 2.9]{DLM1}. (2) is proved in \cite[Proposition 5.4]{DLM1}.
For (3) see \cite[Theorem 3.7]{DLM2}, \cite[Theorem 15.2]{bonnafe action}
and \cite[Theorem 14.11]{bonnafe asterisque}.
\end{proof}

It is conjectured that the above theorem holds without any restriction
(on $p$, $q$, $F$ or the centre of $\Lb$...). However, at the time  of  the writing of this paper, this general conjecture is still unproved.
So we will denote by $\GC(\Gb,\Lb,\Pb)$ the property
$$\lexp{*}{R}_{\Lb \subset \Pb}^\Gb [K\G^\Gb]=(-1)^{d_\Pb} [K\G^\Lb].
\leqno{(\GC(\Gb,\Lb,\Pb))}$$
Most of the results of this subsection will be valid only
under the hypothesis that $\GC(\Gb,\Lb,\Pb)$ holds.  In light of the  above
theorem, and as there are many other  indications  that $\GC(\Gb,\Lb,\Pb)$ holds
in general, this should not be viewed as a big restriction.

In fact, there is also strong evidence that the perfect complex of $\OC\Lb^F$-modules
$\lexp{*}{\RC_{\Lb \subset \Pb}^\Gb} \G^\Gb$ is concentrated in degree $d_\Pb$:
more precisely, it is  conjectured \cite[Conjecture 2.3]{BR2}  that
$$\lexp{*}{\RC}_{\Lb \subset \Pb}^\Gb \G^\Gb \simeq \G^\Lb[-d_\Pb]
\leqno{(\GC^+(\Gb,\Lb,\Pb))}$$
The conjectural property $\GC^+(\Gb,\Lb,\Pb)$ is a far reaching extension of
$\GC(\Gb,\Lb,\Pb)$. It is known to hold only if $\Pb$ is $F$-stable (see Theorem
\ref{rappel restriction} (1) or \cite[Theorem 2.1]{BR2} for a
module-theoretic proof)
or if $\Lb$ is a maximal torus and $(\Pb,F(\Pb))$ lies in the orbit associated with
an element of the Weyl group which is a product of simple reflections
lying in different $F$-orbits \cite[Theorem 3.10]{BR2}. Of course, a proof of this
conjecture would produce immediately a morphism of $\OC$-algebras $\HC^\Gb \to \HC^\Lb$
(which is uniquely determined since $\HC^\Lb$ is commutative). However,
as we shall see in this section, we only need that $\GC(\Gb,\Lb,\Pb)$
holds to get the following result:

\bigskip

\begin{prop}\label{formule trace}
If $\GC(\Gb,\Lb,\Pb)$ holds, then, for all
$h \in K\HC^\Gb \subset K\Gb^F$,
$$\KCur_{\Lb}^\Gb(h) =
(-1)^{d_\Pb} \sum_{s \in \Lb_\sem^{*F^*}/\sim_{\Lb^{*F^*}}}
\Tr_{\Lb \subset \Pb}^\Gb(h,e_{\chi_s^\Lb}e_{\psi_\Lb}) e_{\chi_s^\Lb}e_{\psi_\Lb}.$$
\end{prop}

\begin{proof}
We assume throughout this proof that $\GC(\Gb,\Lb,\Pb)$ holds.
We denote by $\Gamb^\Gb$ the character afforded by the module $K\G^\Gb$.
Let $f : K\HC^\Gb \to K\HC^\Lb$ be the map defined by
$$f(h)=(-1)^{d_\Pb} \sum_{s \in \Lb_\sem^{*F^*}/\sim_{\Lb^{*F^*}}}
\Tr_{\Lb \subset \Pb}^\Gb(h,e_{\chi_s^\Lb}e_{\psi_\Lb}) e_{\chi_s^\Lb}e_{\psi_\Lb}.$$
Let $s \in \Lb_\sem^{*F^*}$. In order to prove the proposition, we only need
to check that
$$\chi_s^{\HC^\Lb} \circ f = \chi_s^{\HC^\Lb} \circ \KCur_\Lb^\Gb.\leqno{(?)}$$

First, note that $\langle \chi_s^\Lb, \G^\Lb \rangle_{\Lb^F}=1$ so,
by adjunction, and since $\GC(\Gb,\Lb,\Pb)$ holds, we have
$\langle R_{\Lb \subset \Pb}^\Gb \chi_s^\Lb, \G^\Gb \rangle_{\Gb^F}=(-1)^{d_\Pb}$.
Since $\chi_s^\Gb$ is the unique irreducible constituent of $\Gamb^\Gb$ lying
in $\EC(\Gb^F,(s)_{\G^{*F^*}})$ and since all the irreducible constituents
of $R_{\Lb \subset \Pb}^\Gb \chi_s^\Lb$ belong to $\EC(\Gb^F,(s)_{\Gb^{*F^*}})$
(see for instance \cite[Theorem 11.10]{bonnafe asterisque}),
we have
$$R_{\Lb \subset \Pb}^\Gb \chi_s^\Lb = (-1)^{d_\Pb} \chi_s^\Gb +
\sum_{\substack{\chi \in \EC(\Gb^F,(s)_{\Gb^{*F^*}}) \\ \chi \neq \chi_s^\Gb}}
m_\chi \chi$$
for some $m_\chi \in \ZM$. By Proposition \ref{representations H} (f), we have
$$\bigl(R_{\Lb \subset \Pb}^\Gb \chi_s^\Lb\bigr)(h)= (-1)^{d_\Pb} \chi_s^\Gb(h)
= (-1)^{d_\Pb} \chi_s^{\HC^\Gb}(h)= \chi_s^{\HC^\Lb}( \KCur_\Lb^\Gb(h))  \leqno{(*)}$$
for all $h \in K\HC^\Gb$. On the other hand, we have
$$\chi_s^{\HC^\Lb}(f(h)) = (-1)^{d_\Pb}
\Tr_{\Lb \subset \Pb}^\Gb(h,e_{\chi_s^\Lb} e_{\psi_\Lb}).$$
But, since the actions of $h$ and of $e_{\chi_s^\Lb} e_{\psi_\Lb}$ on the
cohomology groups $H^i_c(\Yb_\Pb,K)$ commute and since $e_{\chi_s^\Lb} e_{\psi_\Lb}$
is an idempotent, we have that $\Tr_{\Lb \subset \Pb}^\Gb(h,e_{\chi_s^\Lb} e_{\psi_\Lb})$
is the trace of $h$ on the virtual module
$$\sum_{i \ge 0} (-1)^i [H^i_c(\Yb_\Pb,K) e_{\chi_s^\Lb} e_{\psi_\Lb})]
= \sum_{i \ge 0} (-1)^i [H^i_c(\Yb_\Pb,K) \otimes_{K\Lb^F} K\Lb^F
e_{\chi_s^\Lb} e_{\psi_\Lb})].$$
Now, by Proposition \ref{representations H} (e), the $K\Lb^F$-module
$K\Lb^F e_{\chi_s^\Lb} e_{\psi_\Lb}$ affords the character $\chi_s^\Lb$.
So it follows that
$$\Tr_{\Lb \subset \Pb}^\Gb(h,e_{\chi_s^\Lb} e_{\psi_\Lb}) =
\bigl(R_{\Lb \subset \Pb}^\Gb \chi_s^\Lb\bigr)(h).\leqno{(**)}$$
So, $(?)$ follows from the comparison of $(*)$ and $(**)$.
\end{proof}

\begin{prop}\label{curtis tore}
If $\Bb$ is a Borel subgroup of $\Gb$ and if $\Tb$ is
a maximal torus of $\Bb$, then $\GC(\Gb,\Tb,\Bb)$ holds and
$\KCur_\Tb^\Gb$
coincides with Curtis homomorphism $f_\Tb$ defined in \cite[Theorem 4.2]{curtis}.
We have, for all $h \in K\HC^\Gb$,
$$\KCur_\Tb^\Gb (h)=\frac{1}{|\Tb^F|} \sum_{t \in \Tb^F}
\Tr_{\Tb \subset \Bb}^\Gb(h,t) t^{-1}.$$
\end{prop}

\bigskip

\noindent{\sc Remark - } The formula given in Proposition \ref{curtis tore}
gives a concise form for Curtis homomorphism. It can be checked directly,
using the character formula \cite[Proposition 12.2]{dmbook}, that this indeed coincides
with the formulas given by Curtis in terms of Green functions
\cite[4.3]{curtis}. However, we shall give a more theoretical proof of this
coincidence.\finl

\begin{proof}
Since the centre of $\Tb$ is connected, $\GC(\Gb,\Tb,\Bb)$ holds
by Theorem \ref{rappel restriction} (2).
Also, $\Ub_\Tb=1$, $\psi_\Tb=1$, so $K\HC^\Tb = K\Tb^F$. So the primitive
idempotents of $K\HC^\Tb$ are the primitive idempotents of $K\Tb^F$ and
the formula given above can be obtained by a straightforward computation.

Now, let $\Tb^*$ be an $F^*$-stable
maximal torus of $\Gb^*$ dual to $\Tb$. If $s \in \Tb^{*F^*}$, then
$\chi_s^\Tb = \chi_s^{\HC^\Tb}$ and,
by \cite[Theorem 4.2]{curtis}, Curtis homomorphism
$f_\Tb : K\HC^\Gb \to K\Tb^F$ satisfies
$$\chi_s^\Tb \circ f_\Tb = \chi_s^{\HC^\Gb}.$$
Since $\chi^{\HC^\Tb}$ is an isomorphism of $K$-algebras,
we get from Proposition \ref{prop curtis} that $f_\Tb=\Cur_{\Tb}^\Gb$.
\end{proof}

\bigskip

\remark{rlg chi}
If $\chi$ is a class function on $\Lb^F$ (which can be seen as a
class function on $K\Lb^F$) and if $\GC(\Gb,\Lb,\Pb)$ holds, then we have
$$\chi(\Cur_{\Lb}^\Gb(h))=(-1)^{d_\Pb} R_{\Lb \subset \Pb}^\Gb(\chi)(h)$$
for all $h \in K\HC^\Gb$. For this, one may assume that $\chi \in \Irr \Lb^F$.
Let $s \in \Lb_\sem^{*F^*}$ be such that $\chi \in \EC(\Lb^F,(s)_{\Lb^{*F^*}})$.
If $\chi=\chi_s^\Lb$, then this is the equality $(*)$ in the proof of
the Proposition \ref{formule trace}. If $\chi \neq \chi_s^\Lb$, we must
show that $R_{\Lb \subset \Pb}^\Gb(\chi)(h)=0$ for all $h \in K\HC^\Gb$
(see Proposition \ref{representations H} (f)).
Let $\g \in \Irr \Gb^F$ be such that
$\langle \g,R_{\Lb \subset \Pb}^\Gb \chi \rangle_{\Gb^F} \neq 0$.
Then $\g \in \EC(\Gb^F(s)_{\Gb^{*F^*}})$
(see for instance \cite[Theorem 11.10]{bonnafe asterisque}) and, by
Proposition \ref{representations H} (f), it is sufficient to show that
$\g \neq \chi_s^\Gb$. But
$$\langle \chi_s^\Gb,R_{\Lb \subset \Pb}^\Gb \chi \rangle_{\Gb^F} =
\langle \G^\Gb,R_{\Lb \subset \Pb}^\Gb \chi \rangle_{\Gb^F}=
\langle \G^\Lb,\chi  \rangle_{\Lb^F} = 0.$$
This shows the result.\finl

\bigskip

\subsection{A morphism ${\boldsymbol{\HC^\Gb \to \HC^\Lb}}$\label{sub OO}}
We conjecture that, in general, $\KCur_{\Lb}^\Gb(\HC^\Gb) \subset \HC^\Lb$.
At this stage of the paper, we are only able to prove it in the following
cases (in Corollary \ref{cur O}, we shall see that this property
also holds if $\ell$ does not divide the order of $W$):

\bigskip

\begin{theo}\label{main curtis}
We have:
\begin{itemize}
\itemth{a} If $\GC^+(\Gb,\Lb,\Pb)$ holds,
then $\KCur_{\Lb}^\Gb(\HC^\Gb) \subset \HC^\Lb$ and the resulting
morphism of $\OC$-algebra $\HC^\Gb \to \HC^\Lb$ coincides with the functorial morphism
coming from the isomorphism
$\lexp{*}{\RC}_{\Lb \subset \Pb}^\Gb \G^\Gb \simeq \G^\Lb[-d_\Pb]$.

\itemth{b} If $\Lb$ is a maximal torus, then $\KCur_{\Lb}^\Gb(\HC^\Gb)
\subset \HC^\Lb$.
\end{itemize}
\end{theo}

\begin{proof}
(b) follows easily from Proposition \ref{curtis tore} and from the well-known
fact that, if $(g,l) \in \Gb^F \times \Lb^F$, then $|\Lb^F|$ divides
$\Tr_{\Lb \subset \Pb}^\Gb(g,l)$ because $\Lb$ is a maximal
torus.

\medskip

(a)
The complex $\Rrm\G_c(\Yb_\Pb^\Gb,\OC)$ is perfect as a complex of left $\OC\Gb^F$-modules.
Therefore, we have
$\lexp{*}{\RC}_{\Lb \subset \Pb}^\Gb C =  \Rrm\G_c(\Yb_\Pb^\Gb,\OC)^* \otimes_{\OC\Gb^F} C$
for any complex $C$ of $\OC\Gb^F$-modules. If $\GC^+(\Gb,\Lb,\Pb)$ holds, then
this means that we have an isomorphism
$\Rrm\G_c(\Yb_\Pb^\Gb,\OC)^*e_\psi \simeq \G^\Lb[-d_\Pb]$.
In particular, the complex $\Rrm\G_c(\Yb_\Pb^\Gb,\OC)^*e_\psi$ is concentrated
in degree $d_\Pb$. Therefore, there exists an $(\OC\Lb^F,\HC^\Gb)$-bimodule
$P$ such that $\Rrm\G_c(\Yb_\Pb^\Gb,\OC)^*e_\psi \simeq P[-d_\Pb]$.
Moreover, as a left $\OC\Lb^F$-module, we have
an isomorphism $\a : \OC\Lb^F e_{\psi_\Lb} \longmapright{\sim} P$.

This induces a morphism
$$\fonction{\alpt}{\HC^\Gb}{\HC^\Lb}{h}{\a^{-1}(\a(e_{\psi_\Lb}) h).}$$
The morphism $\alpt : \HC^\Gb \to \HC^\Lb$ does not depend on the choice
of $\a$ because $\HC^\Lb$ is commutative. This morphism can be extended to
a morphism $\lind{K}{\alpt} : K\HC^\Gb \to K\HC^\Lb$,
$h \mapsto \lind{K}{\a}^{-1}(\lind{K}{\a}(e_{\psi_\Lb}) h)$.
Now the Theorem would follow
if we show that $\lind{K}{\alpt}=\KCur_{\Lb}^\Gb$.

So let $s \in \Gb_\sem^{*F^*}$. Let $\EC$ be a set of representatives of
$\Lb^F$-conjugacy classes which are contained in $\Lb^{*F^*} \cap (s)_{\Gb^{*F^*}}$
and let $e=\sum_{t \in \EC} e_{\psi_\Lb} e_{\chi_t^\Lb}$.
Then $\lexp{*}{R}_{\Lb \subset \Pb}^\Gb \chi_s^\Gb=(-1)^{d_\Pb} \sum_{t \in \EC}\chi_t^\Lb$.
In particular, $\a$ induces an isomorphism
$$KP e_{\chi_s^\Lb} \simeq K\Lb^F e.$$
So this shows that $\lind{K}{\alpt}(e_\psi e_{\chi_s^\Gb})=e$, as desired.
\end{proof}

\bigskip

\subsection{Truncation at unipotent blocks\label{sub unipotent l}}
We denote by $b^\Gb$ the sum of the {\it unipotent} block idempotents of $\Gb^F$.
In other words,
$$b^\Gb=\sum_{\substack{s \in \Gb^{*F^*}_\sem/\sim \\ \text{$s$ is an $\ell$-element}}}
\sum_{\chi \in \EC(\Gb^F,(s)_{\Gb^{*F^*}})} e_\chi.$$
The algebra $\HC^\Gb$ is a module over the centre of the $\OC$-algebra
$\OC\Gb^F$: so $b^\Gb \HC^\Gb$ is an $\OC$-algebra with unit $b^\Gb e_\psi$.
Note that
$$b^\Gb e_\psi=\sum_{\substack{s \in \Gb^{*F^*}_\sem/\sim \\ \text{$s$ is an $\ell$-element}}}
e_{\chi_s^\Gb}e_\psi.$$
Now, by definition, we get
$$\KCur_\Lb^\Gb(b^\Gb e_\psi) = b^\Lb e_{\psi_\Lb}.$$
In particular,
\equat\label{unip unip}
\KCur_\Lb^\Gb(b^\Gb K\HC^\Gb) \subset b^\Lb K\HC^\Lb.
\endequat
Let us also recall for future reference the following classical fact:

\bigskip

\begin{prop}
The projective $\OC\Gb^F$-module $b^\Gb \G^\Gb$ is indecomposable.
\end{prop}

\begin{proof}   See  \cite[Proposition 19.6 (i)]{cabanes}.  Note that 
the  statement in \cite{cabanes} is made   under the hypotheses   that  $\Gb$  
has connected center,  
but  the proof   applies   without change   in the  general situation.
\end{proof}

\bigskip

\begin{coro}
The algebra $b^\Gb\HC^\Gb$ is local.
\end{coro}

\bigskip

\section{Glueing Curtis homomorphisms for maximal tori}

\medskip

If $\Bb$ is a Borel subgroup of $\Gb$ and if $\Tb$ is an $F$-stable
maximal torus of $\Bb$, we then write $R_\Tb^\Gb$,
and $\Tr_\Tb^\Gb$ for the maps
$R_{\Tb \subset \Bb}^\Gb$ and $\Tr_{\Tb \subset \Bb}^\Gb$
(see the comments at the end of subsection \ref{subsection dl}
and Proposition \ref{prop curtis} (c)).

Let $\Tb_\Gb$ denote an $F$-stable maximal torus of $\Bb_\Gb$. We set
$W=N_\Gb(\Tb_\Gb)/\Tb_\Gb$.
For each $w \in W$, we fix an element $g \in \Gb$ such that
$g^{-1}F(g)$ belongs to $N_\Gb(\Tb_\Gb)$ and represents $w$. We then
set $\Tb_w = g\Tb_\Gb g^{-1}$. We then define, following \cite[Lemma 1]{curtis shoji},
$$\fonction{\Cur^\Gb}{\HC^\Gb}{\DS{\prod_{w \in W}} \OC\Tb_w^F}{h
\vphantom{\DS{\frac{\DS{A^A_A}}{A}}}}{\bigl(\Cur_{\Tb_w}^\Gb(h)\bigr)_{w \in W}.}$$
The aim of this section is to study the map $\Cur^\Gb$.

\bigskip

\subsection{Properties of ${\boldsymbol{\KCur^\Gb}}$}
Before studying $\Cur^\Gb$, we shall study the simpler map $\KCur^\Gb$.
It turns out that $\KCur^\Gb$ is injective and it is relatively easy to describe
its image: both facts were obtained by Curtis and Shoji \cite[Lemmas 1 and 5]{curtis shoji}
but we shall present here a concise proof.

We first need to introduce some notation. If $w \in W$, we fix an $F^*$-stable
maximal torus $\Tb_w^*$ dual to $\Tb_w$. If $t \in \Tb_w^{*F^*}$, then $\chi_t^{\Tb_w}$
is a linear character of $\Tb_w^F$. If $s$ is a semisimple element of
$\Gb^{*F^*}$, we set
$$e^\Gb(s)=
\Bigl(\sum_{t \in (s)_{\Gb^{*F^*}} \cap \Tb_w^{*F^*}} e_{\chi_t^{\Tb_w}} \Bigr)_{w \in W}
~\in \prod_{w \in W} K\Tb_w^F.$$
Then, by definition, we have
\equat\label{image es}
\KCur^\Gb(e_{\chi_s^\Gb}e_\psi)= e^\Gb(s).
\endequat
Since $(e^\Gb(s))_{(s) \in \Gb_\sem^{*F^*}/\sim}$ is a $K$-linearly independent
family in $\prod_{w\in W} K\Tb_w^F$, we get:

\bigskip

\begin{prop}[Curtis-Shoji]\label{injectif image}
The map $\KCur^\Gb$ is injective and
$$\im \KCur^\Gb = \mathop{\oplus}_{(s) \in \Gb_\sem^{*F^*}/\sim} K e^\Gb(s).$$
\end{prop}

\bigskip

\begin{coro}\label{injectif cur}
The map $\Cur^\Gb$ is injective.
\end{coro}

\bigskip

We shall now recall a characterization of elements of the image of $\KCur^\Gb$
which was obtained by Curtis and Shoji \cite[Lemma 5]{curtis shoji}.
We need some notation.
Let $\SC^\Gb$ denote the set of pairs $(w,\th)$ such that $w \in W$ and $\th$ is a linear
character of $\Tb_w^F$ (which may also be viewed as a morphism of algebras
$\OC\Tb_w^F \to \OC$ or $K\Tb_w^F \to K$). If $(w,\th)$ and $(w',\th')$ are
two elements of $\SC^\Gb$, we write $(w,\th) \equiv (w',\th')$ if $(\Tb_w,\th)$ and
$(\Tb_{w'},\th')$ lie in the same rational series (see for instance
\cite[Definition 9.4]{bonnafe asterisque} for a definition).

\begin{coro}[Curtis-Shoji]\label{image curtis}
Let $t=(t_w)_{w \in W} \in \prod_{w \in W} K\Tb_w^F$. Then
$t \in \im \KCur^\Gb$ if and only if, for all $(w,\th)$, $(w',\th') \in \SC^\Gb$ such
that $(w,\th) \equiv (w',\th')$, we have $\th(t_w)=\th'(t_{w'})$.
\end{coro}

\begin{proof} Let $t=(t_w)_{w \in W} \in \prod_{w \in W} K\Tb_w^F$.  
Since  for all $w \in W$,   
$K\Tb_w^F$ is  split commutative  and semi-simple, 
the idempotents of $K\Tb_w ^F$  form a $K$-basis  of $K\Tb_w ^F$, and  
we  may write
$ t = \prod_{w \in W} \sum_{g \in \Tb_w^{*F^*}}  \alpha_g^{\Tb_w} e_{\chi_g^{\Tb}}, $   where $ \alpha_g^{\Tb_w}  \in K  $.     

 Now, from  Proposition 
\ref{injectif image} we have  that  $ t \in \im \KCur^\Gb$ if and only if, whenever $g, g' $ are rationally conjugate  semi-simple elements of $\Gb^{*F^*}$, then  for any $ w, w' \in W$ such that $ g \in \Tb^{*F^*} $ and $g' \in\Tb'^{*F^*}$, we have $ \alpha_g^{\Tb_w} = \alpha_{g'}^{\Tb_{w'}}$.   
On the other hand,  if   $ g \in  \Tb_w^{*F^*}$, then  
$ \alpha_g^{\Tb_w} = \chi_g^{\Tb_w}(t_w)$. The result follows from the 
definition of the equivalence relation on  $\SC^\Gb$.   

\end{proof}

\bigskip

\subsection{Symmetrizing form}
The $\OC$-algebra $\HC^\Gb$ is symmetric. In particular, the $\OC$-algebra
$\im \Cur^\Gb$ is symmetric (see Corollary \ref{injectif cur}).
We shall give in this subsection a precise formula for the
symmetrizing form on $\im \Cur^\Gb$.
For this, we introduce the following symmetrizing form
$$\fonction{\taut}{\DS{\prod_{w \in W}} K\Tb_w^F}{K}{(x_w)_{w \in W}}{\DS{\frac{1}{|W|}
\sum_{w \in W} \t_w(x_w)},}$$
where $\t_w : K\Tb_w^F \to K$ is the canonical symmetrizing form.

We denote by $\t : \OC\Gb^F \to \OC$ the canonical symmetrizing form.
We denote by $\t_\HC$ the restriction of $|\Ub_\Gb^F| \t$ to $\HC^\Gb$: it is
a symmetrizing form on $\HC^\Gb$ (recall that $|\Ub_\Gb^F|$ is
invertible in $\OC$ and is the highest power of $p$ dividing $|\Gb^F|$). Note that
$$\t_\HC(e_\psi)=1.$$
Of course, the extension $\lind{K}{\t_\HC} : K\HC^\Gb \to K$ is a symmetrizing form
on $K\HC^\Gb$. We have
\equat\label{trace}
\lind{K}{\t}_\HC = \taut \circ \KCur^\Gb.
\endequat
\begin{proof}
Since $\t_w$ is a class function on $\Tb_w^F$, we have,
by Remark \ref{rlg chi},
$$\taut(\KCur^\Gb(h))=\frac{1}{|W|} \sum_{w \in W} R_{\Tb_w}^\Gb(\t_w)(h)$$
for all $h \in K\HC^\Gb$. But, by \cite[Proposition 12.9 and Corollary 12.14]{dmbook},
we have
$$\frac{1}{|W|} \sum_{w \in W} R_{\Tb_w}^\Gb(\t_w)= |\Ub_\Gb^F| \t.$$
This completes the proof of the formula \ref{trace}.
\end{proof}

%
%
%
%
%
%
%

\medskip

\subsection{On the image of ${\boldsymbol{\Cur^\Gb}}$}
We are not able to determine in general the sub-$\OC$-algebra
$\im(\Cur^\Gb)$ of $\prod_{w \in W} \OC\Tb_w^F$. Of course,
we have
\equat\label{image inclusion}
\im(\Cur^\Gb) \subset \im(\KCur^\Gb) \cap \bigl(\prod_{w \in W} \OC\Tb_w^F\bigr).
\endequat
However, there are cases where this inclusion is an equality:

\begin{theo}\label{image}
If $\ell$ does not divide the order of $W$, then
$$\im(\Cur^\Gb) =\im(\KCur^\Gb) \cap \bigl(\prod_{w \in W} \OC\Tb_w^F\bigr).$$
\end{theo}

\begin{proof}
Let $A$ be the image of $\Cur^\Gb$. Then, since $\HC^\Gb$ is a symmetric
algebra (with symmetrizing form $\t_\HC$), it follows from \ref{trace} that
$A$ is a symmetric algebra (with symmetrizing form $\taut_A$, the restriction
of $\taut$ to $A$).

Now, let $B=\im(\KCur^\Gb) \cap \bigl(\prod_{w \in W} \OC\Tb_w^F\bigr)$.
If $\ell$ does not divide $|W|$, then the restriction of $\taut$ to $B$
defines a map $\taut_B : B \to \OC$. By construction, we have $A \subset B \subset KA$.
So the result follows from Lemma \ref{symmetric lemma} below.
\end{proof}

\begin{lem}\label{symmetric lemma}
Let $(\AC,\t)$ be a symmetric $\OC$-algebra and let $\BC$ be a subring
of $K\AC$ such that $\AC \subset \BC$ and $\lind{K}{\t}(\BC) \subset \OC$.
Then $\AC=\BC$.
\end{lem}

\begin{proof}
Let $(a_1,\dots,a_n)$ be an $\OC$-basis of $\AC$ and let $(a_1^*,\dots,a_n^*)$
denote the dual $\OC$-basis of $\AC$ (with respect to $\t$). Then, for all
$h \in K\AC$, we have $h=\sum_{i=1}^n \lind{K}{\t}(h a_i^*) a_i$.
Now, if moreover $h \in \BC$, then $ha_i^* \in \BC$ for all $i$, so
$\lind{K}{\t}(h a_i^*) \in \OC$. So $h \in \AC$.
\end{proof}

\remark{ell}
If $\ell$ does not divide the order of $W$, then the Sylow $\ell$-subgroups
of $\Gb^F$ are abelian. If $\Gb^F=\Sb\Lb_2(\FM_{\! q})$, if $q$ is odd and if
$\ell = 2$, then the inclusion \ref{image inclusion} is strict.
If $\Gb^F=\Gb\Lb_3(\FM_2)$ and if $\ell=3$, then
$\ell$ divides $|W|$ but the Sylow $3$-subgroups of $\Gb^F$ are abelian:
in this case, a brute force computation shows that the inclusion
\ref{image inclusion} is an equality. This suggests the following
question: {\it do we have an equality in \ref{image inclusion} if and only
if the Sylow $\ell$-subgroups of $\Gb^F$ are abelian?}

\bigskip

By Corollary \ref{image curtis} and Theorem \ref{image}, we get:

\bigskip

\begin{coro}\label{caracterisation image}
Let $t=(t_w) \in \prod_{w \in W} \OC \Tb_w^F$ and assume that $\ell$ does
not divide the order of $W$. Then $t \in \im \Cur^\Gb$ if and only if,
for all $(w,\th)$, $(w',\th') \in \SC^\Gb$ such
that $(w,\th) \equiv (w',\th')$, we have $\th(t_w)=\th'(t_{w'})$.
\end{coro}

\begin{coro}\label{image dans O}
Let $h \in K\HC^\Gb$ and assume that $\ell$ does not divide the order
of $W$. Then $h \in \HC^\Gb$ if and only if $\KCur_\Tb^\Gb(h) \in \OC\Tb^F$
for all $F$-stable maximal tori of $\Gb$.
\end{coro}

The next result has been announced at the beginning of \S\ref{sub OO}.

\begin{coro}\label{cur O}
If $\Lb$ is an $F$-stable Levi subgroup of a parabolic subgroup of $\Gb$ and
if $\ell$ does not divide the order of $W$, then
$\KCur_\Lb^\Gb(\HC^\Gb) \subset \HC^\Lb$.
\end{coro}

\begin{proof}
 Let $h \in \HC^\Gb$ and let $h'=\KCur_\Lb^\Gb(h)$.
By Corollary \ref{image dans O}, it is sufficient to show
that $\KCur_\Tb^\Lb(h') \in \OC\Tb^F$ for all $F$-stable maximal
torus $\Tb$ of $\Lb$. But this follows from the transitivity
of the Curtis maps (see Proposition \ref{transitivity})
and from the fact that $\KCur_\Tb^\Gb(\HC^\Gb) \subset \OC\Tb^F$
(see Theorem \ref{main curtis}).
\end{proof}

\bigskip

\subsection{Truncation at unipotent blocks\label{sub unipotent t}}
We keep the notation introduced in \S\ref{sub unipotent l}: for instance,
$b^\Gb$ denotes the sum of the unipotent blocks of $\Gb^F$.

\begin{theo}\label{image sylow}
Assume that $\ell$ does not divide the order of $W$. Let $S$ denote
a Sylow $\ell$-subgroup of $\Gb^F$ and let $\Tb$ denote a maximally split
$F$-stable maximal torus of $C_\Gb(S)$. Then $\Cur_\Tb^\Gb$ induces an isomorphism
$$b^\Gb \HC^\Gb \simeq (\OC\Tb^F b^\Tb)^{N_{\Gb^F}(\Tb)} \simeq
(\OC S)^{N_{\Gb^F}(S)}.$$
\end{theo}

\begin{proof}
First, since $\ell$ does not divide the order of $W$,
$S$ is contained in some maximal torus and the centralizer
$C_\Gb(S)$ is an $F$-stable Levi subgroup of a parabolic subgroup of $\Gb$.
In particular, $S$ is abelian, $\Tb^F$ contains $S$ and $S$ is a Sylow $\ell$-subgroup
of $\Tb^F$. This implies that $N_{\Gb^F}(\Tb) \subset N_{\Gb^F}(S)$.
Moreover, if $n \in N_{\Gb^F}(S)$,
then $\lexp{n}{\Tb}$ is another maximally split maximal torus of $C_\Gb(S)$ so there
exists $g \in C_{\Gb^F}(S)$ such that $\lexp{n}{\Tb}=\lexp{g}{\Tb}$. This shows that
$$N_{\Gb^F}(S) = N_{\Gb^F}(\Tb). C_{\Gb^F}(S).\leqno{(1)}$$
This also implies that the map $\OC S \to \OC \Tb^F b^\Tb$, $x \mapsto x b^\Tb$
induces an isomorphism
$$(\OC\Tb^F b^\Tb)^{N_{\Gb^F}(\Tb)} \simeq
(\OC S)^{N_{\Gb^F}(S)}.$$
So we only need to show that $\Cur_\Tb^\Gb$ induces an isomorphism of algebras
$b^\Gb \HC^\Gb \simeq (\OC\Tb^F b^\Tb)^{N_{\Gb^F}(\Tb)}$.

\medskip

Now, by \ref{unip unip}, we have
that $\Cur_\Tb^\Gb(b^\Gb \HC^\Gb) \subset (\OC\Tb^F b^\Tb)^{N_{\Gb^F}(\Tb)}$.
So it remains to prove that $\Cur_\Tb^\Gb$ is injective on $b^\Gb \HC^\Gb$
and that the above inclusion
is in fact an equality.

\medskip

Let us first prove that $\Cur_\Tb^\Gb$ is injective on $b^\Gb \HC^\Gb$.
Let $\Tb^*$ denote an $F^*$-stable maximal torus which is dual to $\Tb$.
Let $S^*$ denote the Sylow $\ell$-subgroup of $\Tb^{*F^*}$.
Then $|\Gb^F|=|\Gb^{*F^*}|$ and $|\Tb^F|=|\Tb^{*F^*}|$ so $S^*$ is a
Sylow $\ell$-subgroup of $\Gb^{*F^*}$. In particular, every $\ell$-element
of $\Gb^{*F^*}$ is conjugate to an element of $S^*$. So
$\KCur_\Tb^\Gb$ is injective on $b^\Gb K\HC^\Gb$, as desired.

Moreover, since $S^*$ is abelian, two elements of $S^*$ are conjugate in
$\Gb^{*F^*}$ if and only if they are conjugate under $N_{\Gb^{*F^*}}(S^*)$ that is,
if and only if they are conjugate under $N_{\Gb^{*F^*}}(\Tb^*)$:
indeed, by the same argument used above for proving $(1)$, we have
$$N_{\Gb^{*F^*}}(S^*)= N_{\Gb^{*F^*}}(\Tb^*).C_{\Gb^*}(S^*).\leqno{(1^*)}$$
In particular,
$$\KCur_\Tb^\Gb(b^\Gb K\HC^\Gb) = (K\Tb^F b^\Tb)^{N_{\Gb^F}(\Tb)}.\leqno{(2)}$$

So, by (2), we only need to prove that,
$$\text{\it if $h \in b^\Gb K\HC^\Gb$ is such that $\KCur_\Tb^\Gb(h) \in \OC\Tb^F$,
then $h \in b^\Gb \HC^\Gb$.}\leqno{(?)}$$
We shall prove (2) by induction on $\dim \Gb$, the case where $\dim \Gb=\dim \Tb$
being trivial. So let
$h \in b^\Gb K\HC^\Gb$ be such that $\KCur_\Tb^\Gb(h) \in \OC\Tb^F$.
Let $w \in W$. By Corollary \ref{image dans O}, we only need to show that
$\KCur_{\Tb_w}^\Gb(h) \in \OC\Tb_w^F$.

Let $S_w$ denote the Sylow $\ell$-subgroup
of $\Tb_w^F$. Since $S$ is a Sylow $\ell$-subgroup of $\Gb^F$,
we may, and we will, assume that $S_w \subset S$.
Now, let $\Lb=C_\Gb(S_w)$. Since $S_w$ is an $\ell$-subgroup
and $\ell$ does not divide the order of $W$, $\Lb$ is an $F$-stable
Levi subgroup of a parabolic subgroup of $\Gb$. Moreover, we have
$\Tb \subset \Lb$ and $\Tb_w \subset \Lb$.

Now, let $h'=\KCur_\Lb^\Gb(h)$. Then $h' \in b^\Lb K\HC^\Lb$
(see \ref{unip unip}) and, by hypothesis, we have
$\KCur_\Tb^\Lb(h')=\KCur_\Tb^\Gb(h) \in \OC\Tb^F$.
So, if $\dim \Lb < \dim \Gb$, then
$h' \in b^\Lb\HC^\Lb$ by induction hypothesis,
and so $\KCur_{\Tb_w}^\Gb(h)=\KCur_{\Tb_w}^\Lb(h') \in \OC\Tb_w^F$, as
desired. This means that we may, and we will, assume that $\Lb=\Gb$
(or, in other words, that $S_w$ is central in $\Gb$). This implies in particular
that $S_w$ is the Sylow $\ell$-subgroup of $\Zb(\Gb)^F$. Moreover,
since $\ell$ does not divide $|W|$, it does not divide $|\Zb(\Gb)/\Zb(\Gb)^\circ|$,
so $S_w$ is the Sylow $\ell$-subgroup of $(\Zb(\Gb)^\circ)^F$. Since
$|\Tb_w^F|=|\Tb_w^{*F^*}|$ and $|(\Zb(\Gb)^\circ)^F|=|(\Zb(\Gb^*)^\circ)^{F^*}|$,
the Sylow $\ell$-subgroup of $\Tb_w^{*F^*}$ (which we shall denote
by $S_w^*$) is central in $\Gb^*$.

So, let us write
$$h=\sum_{\stackrel{(s) \in \Gb_\sem^{*F^*}/\sim}{\text{$s$ is an $\ell$-element}}} a_s e_{\chi^\Gb_s}e_\psi.$$
Then, by hypothesis,
$$\sum_{s \in S^*} a_s e_{\chi_s^\Tb} \in \OC\Tb^F.$$
In other words, we have, for all $t \in \Tb^F$,
$$\frac{1}{|S|} \sum_{s \in S^*} a_s \chi_s^\Tb(t) \in \OC.\leqno{(3)}$$
We want to show that, for all $t \in \Tb_w^F$,
$$\frac{1}{|S_w|} \sum_{s \in S_w^*} a_s \chi_s^{\Tb_w}(t) \in \OC.\leqno{(??)}$$
Since $\chi_s^{\Tb_w}(t)=1$ if $t$ is an $\ell'$-element of $\Tb_w^F$
and $s \in S_w^*$, we only need to show (??) whenever $t \in S_w$. But,
in this case, $\chi_s^{\Tb_w}(t)=\chi_s^\Tb(t)$ since $t$ is central
in $\Gb$. On the other hand,
let $S'=\{t' \in S~|~\forall~s \in S_w^*,~\chi_s^\Tb(t')=1\}$.
Then $S = S' \times S_w$. So, by (3), we have, forall $t \in S_w$,
$$\frac{1}{|S|}\sum_{t' \in S'}\Bigl(\sum_{s \in S^*} a_s \chi_s^\Tb(tt')\Bigr) \in \OC.$$
But
$$\frac{1}{|S|}\sum_{t' \in S'}\Bigl(\sum_{s \in S^*} a_s \chi_s^\Tb(tt')\Bigr)
=\frac{1}{|S_w|} \sum_{s \in S_w^*} a_s \chi_s^{\Tb_w}(t),$$
so (??) follows.
\end{proof}

\bigskip

\section{The Curtis-Shoji homomorphism}

\medskip

Let $d $ be a fixed positive integer.
In \cite[Theorem 1]{curtis shoji}, Curtis and Shoji defined an
algebra homomorphism from the endomorphism ring of a  Gelfand-Graev
representation of  $K{\Gb^{F^d}}$ to the endomorphism ring of a  Gelfand-Graev
representation of  $K{\Gb^{F}}$.     In this section, we   review the
definition   of this   homomorphism.  We conjecture that  this homomorphism is
defined over $\OC$ and prove this in a special case.

Since we are working with two different isogenies $F$ and $F^d$, we shall
need to use more precise notation. We shall use the index $?_{(e)}$ to denote
the object $?$ considered with respect to the isogeny $F^e$: for instance,
$\G^\Gb_{(d)}$ shall denote a Gelfand-Graev representation of $\Gb^{F^d}$, $\chi_{s,(1)}^\Lb$ shall denote the character $\chi_s^\Lb$ of the finite
group $\Lb^F$ and so on.

\subsection{Notation}
According to our convention,
the regular linear character $\psi$ of $\Ub^F$
will be denoted by $\psi_{(1)}$. We fix a  regular linear character
$\psi_{(d)} : \Ub_\Gb^{F^d} \to \OC^\times \subset K^\times$.
Set
$$\G^\Gb_{(d)} = \OC\Gb^{F^d} e_{\psi_{(d)}} \simeq
\Ind_{\Ub_\Gb^{F^d}}^{\Gb^{F^d}} \OC_{\psi_{(d)}}.$$
and  let  $\HC_{(d)}^\Gb$ denote the endomorphism algebra
of the $\OC\Gb^{F^d}$-module $\G_{(d)}^\Gb$.
For $t \in \Gb_\sem^{*{F^{*d}}}$,
we denote by $\chi_{t,(d)}^\Gb$ the unique
element of $\EC(\Gb^{F^d},(t)_{\Gb^{*{F^{*d}}}})$ which is an irreducible
component of the character afforded by $K\G_{(d)}^\Gb$.
If $\Tb$ is an $F^d$-stable maximal torus, we shall denote
by $\Cur_{\Tb,(d)}^\Gb : \HC^\Gb_{(d)} \to \OC\Tb^{F^d}$
the Curtis homomorphism.

\bigskip

\noindent{\sc Remark - }
By Remark \ref{several}, the endomorphism algebra
$\HC_{(d)}^\Gb$ does not depend on the choice of the regular
linear character $\psi_{(d)}$. There is nevertheless a ``natural''
choice for $\psi_{(d)}$, which is compatible with the theory
of Shintani descent. It is defined as follows.
Consider the map $N : \Ub_\Gb^{F^d}/\Db(\Ub_\Gb)^{F^d} \to \Ub_\Gb^F/\Db(\Ub_\Gb)^F$,
$u \mapsto u F(u) \cdots F^{d-1}(u)$. Then one can take
$\psi_{(d)}=\psi_{(1)} \circ N$.\finl

\bigskip
\subsection {The Curtis-Shoji homomorphism}
For an $F$-stable torus $\Tb$ of  $\Gb$,  denote by
$$\Norm_{F^d/F}^\Tb : \Tb^{F^d} \to \Tb^F, $$ the  surjective group
homomorphism
$$  t \rightarrow  t\cdot\lexp{F}{t} \cdots \lexp{F^{d-1}}{t}. $$
Denote by
$\Norm_{F^d/F}^\Tb$  also the $\OC$-linear  map
$\OC\Tb^{F^d} \to \OC\Tb^{F}$
extending $\Norm_{F^d/F,\Tb}$.

\begin{prop}[Curtis-Shoji]\label{prop curtis shoji original}
There exists a homomorphism of algebras
$$  \KCurSh^\Gb : K\HC_{(d)}^\Gb  \to  K\HC_{(1)}^\Gb $$   which  is  characterized
as the unique linear map from $K\HC_{(d)}^\Gb$ to  $K\HC_{(1)}^\Gb$ with the
property  that
$$  \KCur_{\Tb,(1)}^\Gb \circ  \KCurSh^\Gb =
\KNorm_{F^d/F}^\Tb \circ  \KCur_{\Tb,(d)}^\Gb,  $$
for any  $F$-stable torus   $\Tb$  of  $\Gb$.
\end{prop}

\begin{proof}
See \cite[Theorem 1]{curtis shoji}: the proof uses essentially only the
fact that the map $\KCur_{(1)}^\Gb$ is injective (see Proposition \ref{injectif image})
and the computation of its image (see Corollary \ref{image curtis}).
\end{proof}

\begin{coro}\label{commutation levi}
Let $\Lb$ be an $F$-stable Levi subgroup of a parabolic subgroup of $\Gb$. Then $\KCurSh^\Lb \circ \KCur_{\Lb,(d)}^\Gb = \KCur_{\Lb,(1)}^\Gb \circ \KCurSh^\Lb$.
In other words, the diagram
$$\diagram
\HC_{(d)}^\Gb \rrto^{\DS{\KCur_{\Lb,(d)}^\Gb}} \ddto_{\DS{\KCurSh^\Gb}} &&
\HC_{(d)}^\Lb \ddto^{\DS{\KCurSh^\Lb}} \\
&& \\
\HC_{(1)}^\Gb \rrto^{\DS{\KCur_{\Lb,(1)}^\Gb}} &&\HC_{(1)}^\Lb
\enddiagram$$
is commutative.
\end{coro}

\begin{proof}
Let $f = \KCurSh^\Lb \circ \KCur_{\Lb,(d)}^\Gb$ and
$g = \KCur_{\Lb,(1)}^\Gb \circ \KCurSh^\Lb$.
By Proposition \ref{injectif image}, it is sufficient to
show that $\KCur_{\Tb,(1)}^\Lb \circ f = \KCur_{\Tb,(1)}^\Lb \circ g$
for any $F$-stable maximal torus $\Tb$ of $\Lb$. But this follows
from the transitivity of the Curtis homomorphisms (see Proposition
\ref{transitivity}) and the defining property of the homomorphisms
$\KCurSh^?$ (see Proposition \ref{prop curtis shoji original}).
\end{proof}

We also derive a concrete formula for the map $\KCurSh^\Gb$:

\begin{coro}\label{shintani formule}
The map $\KCurSh^\Gb$ is the unique linear map from
$K\HC_{(d)}^\Gb$ to  $K\HC_{(1)}^\Gb$ with the
property  that  for any  $t \in \Gb_\sem^{*{F^{*d}}}$,
$$ \KCurSh^\Gb(e_{{\chi}_{t,(d)}^\Gb} e_{\psi_{(d)}})=
\sum_{\substack{(s)_{\Gb^{*F^*}} \in \Gb_\sem^{*F^*}/\sim_{\Gb^{*F^*}} \\
(s)_{\Gb^{*F^*}} \subset (t)_{\Gb^{*F^{*d}}}} }e_{{\chi}_{s,(1)}^\Gb} e_{\psi_{(1)}}, $$
In particular,  if
$(t)_{\Gb^{*{F^{*d}}}} \cap \Gb^{*{F^*}}$  is empty,  then
$ \KCurSh( e_{{\chi'}_s^\Gb} e_{\psi'})=0 $.
\end{coro}

\begin{proof}
Let $a= \KCurSh^\Gb(e_{{\chi}_{t,(d)}^\Gb} e_{\psi_{(d)}})$ and
$$b=\sum_{\substack{(s)_{\Gb^{*F^*}} \in \Gb_\sem^{*F^*}/\sim_{\Gb^{*F^*}} \\
(s)_{\Gb^{*F^*}} \subset (t)_{\Gb^{*F^{*d}}}} }e_{{\chi}_{s,(1)}^\Gb} e_{\psi_{(1)}}.$$
By Proposition \ref{injectif image}, we only need to show that,
if $\Tb$ is an $F$-stable maximal torus of $\Gb$, then
$\KCur_{\Tb,(1)}^\Gb(a)=\KCur_{\Tb,(1)}^\Gb(b)$.
But, by Proposition \ref{prop curtis shoji original},
we have
$$\KCur_{\Tb,(1)}^\Gb(a)=\KNorm_{F^d/F}^\Tb(\KCur_{\Tb,(d)}^\Gb
(e_{{\chi}_{t,(d)}^\Gb} e_{\psi_{(d)}}).$$
Therefore,
$$\KCur_{\Tb,(1)}^\Gb(a)=\KNorm_{F^d/F}^\Tb\Bigl(
\sum_{s \in \Tb^{F^{*d}}\cap (t)_{\Gb^{*F^{*d}}}} e_{\chi_{t,(d)}^\Tb}\Bigr).$$
On the other hand,
$$\KCur_{\Tb,(1)}^\Gb(b)=\sum_{s \in \Tb^{F^{*}}\cap (t)_{\Gb^{*F^{*d}}}} e_{\chi_{t,(1)}^\Tb}.$$
So it remains to show that, if $s \in \Tb^{F^{*d}}$,
then
$$\KNorm_{F^d/F}^\Tb(e_{\chi_{s,(d)}^\Tb})=
\begin{cases}
e_{\chi_{s,(1)}^\Tb} & \text{if $s \in \Tb^{*F^*}$,} \\
0 & \text{otherwise.}
\end{cases}$$
But this follows easily from the fact that, by definition \cite[5.21.5, 5.21.6]{delu},
we have $\chi_{s,(1)}^\Tb \circ \Norm_{F^d/F}^\Tb = \chi_{s,(d)}^\Tb$ as
linear characters of $\Tb^{F^d}$.
\end{proof}

\bigskip

\subsection{ A map ${\boldsymbol{\HC_{(d)}^\Gb  \to  \HC_{(1)}^\Gb}}$}
We conjecture that, in general, $\KCurSh^\Gb(\HC_{(d)}^\Gb)\subseteq\HC_{(1)}^\Gb$,
so that $\KCurSh^\Gb$ is defined over $\OC$. However,
we are only able to prove this  in the following special case.

\begin{theo}\label{curtis shoji over dvr}
If $\ell$ does not divide $|W|$,
then  $\KCurSh^\Gb(\HC_{(d)}^\Gb ) \subseteq  \HC_{(1)}^\Gb $.
 \end{theo}

\begin{proof}
Let $a \in \HC_{(d)}^\Gb$ and let $h=\KCurSh^\Gb(a) \in K\HC_{(1)}^\Gb$.
By Corollary \ref{image dans O}, we need to show that, if $\Tb$
is an $F$-stable maximal torus of $\Gb$, then $\KCur_{\Tb,(1)}^\Gb(h) \in \OC\Tb^F$.
But, by Proposition \ref{prop curtis shoji original},  $\KCur_{\Tb,(1)}^\Gb(h)=\KNorm_{F^d/F}^\Tb(t)$, where
$t = \KCur_{\Tb,(d)}^\Gb(a)$. Now, by Theorem \ref{image},
$\KCur_{\Tb,(d)}^\Gb(a) \in \OC\Tb^{F^d}$. So the result
follows from the fact that $\KNorm_{F^d/F}^\Tb$ is defined over $\OC$.
\end{proof}

\bigskip

\subsection{Truncating at unipotent blocks}
Here we  study the    restriction   of   the  Curtis-Shoji homomorphism
to  the  component $b^\Gb_{(d)}K\HC^\Gb$ of $K\HC^\Gb$
(by our usual convention in this section, $b^\Gb_{(m)}$ denotes the sum of the unipotent
block idempotents of $\Gb^{F^m}$).

It is immediate from Corollary \ref{shintani formule}  that
$\KCurSh(b_{(d)}^\Gb e_{\psi_{(d)}}) = b_{(1)}^\Gb e_{\psi_{(1)}}$.
We denote by
$$\KCurSh_\ell^\Gb:  b^\Gb_{(d)}K\HC^\Gb_{(d)} \to b^\Gb_{(1)}K\HC^\Gb_{(1)}, $$
the map obtained by restricting $\KCurSh^\Gb$.

\begin{prop} \label{prop injective surjective}
We have
\begin{itemize}
\itemth{a}  $\KCurSh_\ell^\Gb$ is
surjective if and only if  whenever a   pair  of  $\ell$-elements  of
$\Gb^{*^{F^*}}$ are conjugate in   $\Gb^{*{F^{*d}}}$, they  are also
conjugate in $\Gb^{*{F^{*}}} $.

\itemth{b} $\KCurSh_\ell^\Gb$ is injective if and only if  every  $\ell$-element  of
$\Gb^{*{F^{*d}}} $   is $\Gb^{*{F^{*d}}}$-conjugate to an element of
$\Gb^{*{F^{*}}} $.
\end{itemize}
\end{prop}

\begin{proof}
$\KCurSh_\ell^\Gb $ is a unitary  map of  commutative split
semi-simple algebras, hence is
surjective if and only if  the image of any  primitive idempotent  is either
a primitive idempotent or  $0$. Similarly,  $\KCurSh_\ell^\Gb$ is injective if  the
image of  every  idempotent is  non-zero.  Both parts of the proposition
are  now  immediate
from Corollary \ref{shintani formule}.
\end{proof}

Let $\ZC(\Gb)$ denote the finite group $\Zb(\Gb)/\Zb(\Gb)^\circ$.
The following corollary is related to \cite[Lemma 6]{curtis shoji}:

\begin{coro}\label{centre surjectif}
Let $r$ denote the order of the automorphism of $\ZC(\Gb)_\ell$ induced by $F$.
If ${\mathrm{gcd}}(d,r\ell)=1$, then $\KCurSh_\ell^\Gb$ is surjective.
\end{coro}

\begin{proof}
The proof is somewhat similar to the proof of \cite[Lemma 6]{curtis shoji}:
since our situation is a  bit different and since our hypothesis
is slightly weaker, we shall recall a proof.
Let $s$ and $t$ be two $\ell$-elements of $\Gb^{*F^*}$ and assume that
they are conjugate in $\Gb^{*F^{*d}}$. By Proposition \ref{prop injective surjective} (a),
we only need to show that they are conjugate in $\Gb^{*F^*}$.

So let $g \in \Gb^{*F^{*d}}$ be such that $t=gsg^{-1}$.
Let $A=C_{\Gb^*}(s)/C_{\Gb^*}^\circ(s)$ and let $\s$ denote the
automorphism of $A$ induced by $F^*$. We set $\Ati=A \rtimes <\s>$.
It is a classical fact that $A$ is an $\ell$-group (since $s$ is
an $\ell$-element: see for instance \cite[Lemma 2.1]{broue michel})
and that there is an injective morphism $A \injto \ZC(\Gb)^\wedge$
commuting with the actions of the Frobenius endomorphisms
(see for instance \cite[8.4]{bonnafe asterisque}).
In particular, the order of $\s$ divides $r$. So ${\mathrm{gcd}}(d,\Ati)=1$.
Therefore, the map $\Ati \to \Ati$, $x \mapsto x^d$ is bijective.

Now, since $s$ and $t$ are $F^*$-stable, the element $h=g^{-1}F^*(g)$
belongs to $C_{\Gb^*}(s)$. We denote by $x$ its class in $A$.
The fact that $g$ belongs to $\Gb^{*F^{*d}}$ implies that
$h F^*(h) \cdots F^{*d-1}(h)=1$. So $x \s(x) \cdots \s^{d-1}(x)=1$.
In other words, $(x\s)^d=\s^d$. So $x=1$. In other words,
$g^{-1} F^*(g) \in C_{\Gb^*}^\circ(s)$. By Lang's Theorem,
this implies that $s$ and $t$ are conjugate in $\Gb^{*F^*}$.
\end{proof}

\begin{coro}\label{ordre injectif}
If $\ell$ does not divide $[\Gb^{F^d}:\Gb^F]$, then
$\KCurSh_\ell^\Gb$ is injective.
\end{coro}

\begin{proof}
This follows from the fact that $|\Gb^F|=|\Gb^{*F^*}|$ (and similarly
for $F^d$) and from proposition \ref{prop injective surjective}.
\end{proof}

Let us make a brief comment on this last result. If $r$ denotes the order
of the automorphism induced by $F$ on $\ZC(\Gb)_\ell$ (as in Corollary
\ref{centre surjectif}), it is not clear if the condition that
$\ell$ does not divide $[\Gb^{F^d}:\Gb^F]$ implies that $d$ is prime
to $r$. However, one can easily get the following result:

\begin{lem}\label{ell divise}
If $\ell$ divides $|\Gb^F|$ and does not divide $[\Gb^{F^d}:\Gb^F]$, then
$\ell$ does not divide $d$.
\end{lem}

\begin{proof}
It is sufficient to show that, if $\ell$ divides $|\Gb^F|$, then
$\ell$ divides $[\Gb^{F^\ell}:\Gb^F]$.
For this, let $q_0=q^{1/\d}$ (recall that $F^\d$ is
a Frobenius endomorphism on $\Gb$ with respect to some $\gfq$-structure
on $\Gb$). We denote by $\phi$ the automorphism of $V=X(\Tb) \otimes_\ZM K$
such that $F=q_0 \phi$. Then $\phi$ normalizes $W$ so the invariant algebra
$S(V^*)^W$ can be generated by homogeneous polynomials $f_1$,\dots, $f_n$
(where $n=\dim_K V=\dim \Tb$) which are eigenvectors of $\phi$.
Let $d_i$ denote the degree of $f_i$ and let $\e_i \in K^\times$
be such that $\phi(f_i)=\e_i f_i$. Then
$$|\Gb^F|=q_0^{|\Phi_+|} \prod_{i =1}^n (q_0^{d_i}-\e_i)$$
$$|\Gb^{F^\ell}|=q_0^{\ell|\Phi_+|} \prod_{i =1}^n
(q_0^{\ell d_i}-\e_i^\ell).\leqno{\text{and}}$$
In particular, we have
$$[\Gb^{F^\ell}:\Gb^F] = q_0^{(\ell-1)|\Phi_+|}
\prod_{i=1}^n
(q_0^{d_i(\ell-1)} + q_0^{d_i(\ell-2)} \e_i + \cdots + q_0^{d_i} \e_i^{\ell-2} +
\e_i^{\ell -1}).$$
View this last equality in $\OC$ (and recall that $\lG$ denotes
the maximal ideal of $\OC$).
Now, if $\ell$ divides $|\Gb^F|$, there exists $i$ such that
$$q_0^{d_i} \equiv \e_i \mod \lG.$$
Therefore,
$$q_0^{d_i(\ell-1)} + q_0^{d_i(\ell-2)} \e_i + \cdots +
q_0^{d_i} \e_i^{\ell-2} + \e_i^{\ell -1}
\equiv \ell \e_i^{\ell-1} \equiv 0 \mod \lG.$$
This shows that $[\Gb^{F^\ell}:\Gb^F] \in \lG \cap \ZM = \ell \ZM$.
\end{proof}

If $\KCurSh_\ell^\Gb(b_{(d)}^\Gb \HC^\Gb_{(d)}) \subset b^\Gb_{(1)} \HC^\Gb_{(1)}$,
we denote by $\CurSh_\ell^\Gb : b_{(d)}^\Gb \HC^\Gb_{(d)} \to b^\Gb_{(1)} \HC^\Gb_{(1)}$
the induced map. This happens for instance if $\ell$ does not divide $|W|$
(see Theorem \ref{curtis shoji over dvr}).

\begin{theo}\label{endo iso}
If $\ell$ does not divide $|W|\cdot [\Gb^{F^d}:\Gb^F]$, then
$\CurSh_\ell^\Gb : b_{(d)}^\Gb \HC^\Gb_{(d)} \to b^\Gb_{(1)} \HC^\Gb_{(1)}$
is an isomorphism of algebras.
\end{theo}

\begin{proof}

By Theorem \ref{curtis shoji over dvr}, the map $\CurSh_\ell^\Gb$
is well-defined. By Corollary \ref{ordre injectif}, it is injective.
So it remains to show that it is surjective.

First, the order of $\ZC(\Gb)$ divides the order of $W$.
So, since $\ell$ does not divide the order of $W$, we get that $\ZC(\Gb)_\ell=1$.
So, by Corollary \ref{centre surjectif}, the map $\KCurSh_\ell^\Gb$ is surjective.
So, if $h \in b^\Gb_{(1)} \HC^\Gb_{(1)}$, there exists
$\hti \in b_{(d)}^\Gb K\HC^\Gb_{(d)}$ such that $\KCurSh_\ell^\Gb(\hti)=h$.
So it remains to show that $\hti \in b_{(d)}^\Gb \HC^\Gb_{(d)}$.

Let $S$ be a Sylow $\ell$-subgroup of $\Gb^F$. By hypothesis, it is
a Sylow $\ell$-subgroup of $\Gb^{F^d}$.
Let $\Tb$ be a maximally split $F$-stable maximal torus of $C_\Gb(S)$
(as in Theorem \ref{image sylow}). Let $\tti = \KCur_{\Tb,(d)}^\Gb(\hti)$
and $t=\KCur_{\Tb,(1)}^\Gb(h)$. Then, by Proposition \ref{prop curtis shoji original},
we have
$$t=\KNorm_{F^d/F}^\Tb(\tti)$$
and, by \ref{unip unip},
$$\tti \in K\Tb^{F^d} b^\Tb_{(d)}\qquad\text{and}\qquad t \in K\Tb^F b_{(1)}^\Tb.$$
Also, by the statement $(?)$ of the proof of Theorem \ref{image sylow},
it is sufficient to show that $\tti \in \OC\Tb^{F^d} b^\Tb_{(d)}$.

Write $\tti=\sum_{\zti \in \Tb^{F^d}} a_\zti \zti$ and $t=\sum_{z \in \Tb^F} b_z z$
with $a_\zti \in K$ and $b_z \in \OC$.
Let $H$ be the kernel of the group homomorphism $\Norm_{F^d/F}^\Tb$.
By hypothesis, $S$ is also a Sylow $\ell$-subgroup of $\Tb^{F^d}$.
So $\ell$ does not divide $[\Tb^{F^d}:\Tb^F]=|H|$.
Now, if $\zti \in \Tb^{F^d}$, and if we set $z=\Norm_{F^d/F}^\Tb(\zti)$,
then $b_z =\sum_{h \in H} a_{h\zti}$. But, since $\tti \in K\Tb^{F^d} b^\Tb_{(d)}$,
we have $a_{h\zti}=a_\zti$ for every $h \in H$ (in fact, for every $\ell'$-element
$h$ of $\Tb^{F^d}$). So $|H| a_\zti = b_z \in \OC$,  which means that  $a_\zti \in \OC$
since $|H|$ is invertible in $\OC$.
\end{proof}

\noindent{\sc Acknowledgements.}  The authors would like to thank Raphael Rouquier 
for  useful discussions  on the subject of this paper and especially for 
highlighting the importance of  the  existence of symmetrizing forms.  
They also thank Marc Cabanes and the referee for  their comments and  suggestions.

}

\end{document}